\documentclass{elsart}

\usepackage{amsmath,amssymb, mathrsfs}
\usepackage[T1]{fontenc}
\usepackage{yfonts}

\begin{document}

\begin{frontmatter}
\hsize=6.5in


\title{Identities for the Hankel transform and their applications}

\author[label1]{Ahmet Dernek},
\author[label1]{Ne\c se Dernek},
\author[label2]{Osman Y\"urekli\corauthref{cor1}}

\address[label1]{Department of Mathematics, University of Marmara, TR-34722, Kadik\"oy, Istanbul, Turkey}
\address[label2]{Department of Mathematics, Ithaca College, Ithaca, NY 14850, US}

\corauth[cor1]{Corresponding author}

\begin{abstract}

In the present paper the authors show that  iterations of the Hankel transform with $\mathscr{K}_{\nu}$-transform  is a constant multiple of the Widder transform. Using these iteration identities,  several Parseval-Goldstein type theorems for these transforms are given. By making use of these results  a number of new Goldstein type exchange identities are obtained for these and the Laplace transform. The identities proven in this paper are shown to give rise to useful corollaries for evaluating infinite integrals of special functions. Some examples are also given as illustration of the results presented here.

\end{abstract}


\end{frontmatter}

\section{Introduction, definition and preliminaries}
Over four decades ago, Widder \cite{Wd66} presented a systematic account of the so-called Widder transform:
\begin{equation}
\label{wpt}
\mathscr{P}\big\{f(x);y\big\}
=\int_0^\infty \frac{x\,f(x)}{x^{2}+y^{2}}\,dx,
\end{equation}
which, by an exponential change of variable, becomes a convolution transform with a kernel belonging to a general class treated by Hirschman and Widder \cite{HW}. The Widder transform \eqref{wpt} is related to the classical Stieltjes transform as follows:
\begin{equation}
\label{wt:st}
\mathscr{P}\big\{f(x);y\big\}=\frac{1}{2}\,\mathscr{S}\big\{f\big(\sqrt{x}\big);y^{2}\big\},
\end{equation}
where  the Stieltjes transform is defined by
\begin{equation}
\label{st}
\mathscr{S}\big\{f(x);y\big\}
=\int_0^\infty \frac{f(x)}{x+y}\,dx.
\end{equation}
Srivastava and Singh \cite{SS} gave the following Goldstein type exchange identity for the Widder potential transform \eqref{wpt} together with several illustrative examples of its applications (see also \cite{Go}):
\begin{equation}
\label{eip}
\int_0^\infty y\,\mathscr{P}\big\{f(x);y\big\}\,g(y)\,dy
=\int_0^\infty x\,f(x)\,\mathscr{P}\big\{g(y);x\big\}\,dx.
\end{equation}
The following identity similar to \eqref{eip} was introduced earlier by Goldstein \cite{Go}:
\begin{equation}
\label{eil}
\int_0^\infty \mathscr{L}\big\{f(x);y \big\}\,g(y)\,dy
=\int_0^\infty f(x)\,\mathscr{L} \big\{g(y);x \big\}\,dx 
\end{equation}
is popularly known as the Goldstein type exchange identity for the classical Laplace transform:
\begin{equation}
\label{lt}
\mathscr{L}\big\{f(x);y\big\}
=\int_0^\infty \exp(-x\,y)\,f(x)\,dx.
\end{equation}

Srivastava and Y\"{u}rekli \cite{SY91, SY95} gave various Parseval-Goldstein type identities for the Laplace transform \eqref{lt}, the Widder potential transform \eqref{wpt}, the Fourier sine transform:
\begin{equation}
\label{fst}
\mathscr{F}_{S}\big\{f(x);y \big\}
=\int_0^\infty \sin(x\,y)\,f(x)\,dx, 
\end{equation}
and the Fourier cosine transform:
\begin{equation}
\label{fct}
\mathscr{F}_{C} \big\{f(x);y \big\}
=\int_0^\infty \cos(x\,y)\,f(x)\,dx. 
\end{equation}
Some results involving the Hankel transform and the $\mathscr{K}$-transform were given by (among others) Srivastava and Y\"urekli \cite{SY95}. The Hankel transform is defined by 
\begin{equation}
\label{ht}
\mathscr{H}_{\nu}\big\{f(x);y \big\}
=\int_0^\infty \sqrt{x\,y}\,\text{J}_{\nu}(x\,y)\,f(x)\,dx,
\end{equation}
where $\text{J}_{\nu}(x)$ denotes the Bessel function of the first kind of order $\nu$. Using the formula (cf. \cite[p. 306, Eq. 32:13:10]{SO})
\begin{equation}
\label{hns}
\text{J}_{\frac{1}{2}}(x)=\sqrt{\frac{2}{\pi\,x}}\,\sin(x)
\end{equation}
the definition \eqref{fst} of the Fourier sine transform, and the definition \eqref{ht} of the Hankel transform, we obtain the familiar relationship
\begin{equation}
\label{ht:fst}
\mathscr{H}_{\frac{1}{2}}\big\{f(x);y\big\}=\sqrt{\frac{2}{\pi}}\,\mathscr{F}_{S}\big\{f(x);y\big\}.
\end{equation}
Similarly, using the formula (cf. \cite[p. 306, Eq. 32:13:11]{SO})
\begin{equation}
\label{hnc}
\text{J}_{-\frac{1}{2}}(x)=\sqrt{\frac{2}{\pi\,x}}\,\cos(x)
\end{equation}
the definition \eqref{fct} of the Fourier cosine transform, and the definition \eqref{ht} of the Hankel transform, we obtain the  relationship
\begin{equation}
\label{ht:fct}
\mathscr{H}_{-\frac{1}{2}}\big\{f(x);y\big\}=\sqrt{\frac{2}{\pi}}\,\mathscr{F}_{C}\big\{f(x);y\big\}.
\end{equation}
The $\mathscr{K}$-transform is defined by
\begin{equation}
\label{kt}
\mathscr{K}_{\nu} \big\{f(x);y \big\}
=\int_0^\infty \sqrt{x\,y}\,\text{K}_{\nu}(x\,y)\,f(x)\,dx, 
\end{equation}
where $\text{K}_{\nu}(x)$ is the modified Bessel function (or the Macdonald function) of order $\nu$. Using the formula (cf. \cite[p. 239, Eq. 26:13:5]{SO})
\begin{equation}
\label{knl}
\text{K}_{\frac{1}{2}}(x)=\sqrt{\frac{\pi}{2\,x}}\,\exp(-x),
\end{equation}
the definition \eqref{lt} of the Laplace transform, and the definition \eqref{kt} of the Hankel transform, we obtain the relationship
\begin{equation}
\label{kt:lt}
\mathscr{K}_{\frac{1}{2}}\big\{f(x);y\big\}=\sqrt{\frac{\pi}{2}}\,\mathscr{L}\big\{f(x);y\big\},
\end{equation}
which incidentally holds true also when $\mathscr{K}_{\frac{1}{2}}$ is replaced by $\mathscr{K}_{-\frac{1}{2}}$.

In this article, we first present iteration identities involving the Widder potential transform, the Hankel transform, and $\mathscr{K}$-transform. By using these iteration identities, we establish Parseval-Goldstein type identities involving these transforms. The Parseval-Goldstein type identities established here yields new identities for the various integral transforms introduced above. As applications of the resulting identities and theorems, some illustrative examples are also given.

\section{The main Parseval-Goldstein type theorem}

The following identities involving the Hankel transform, the $\mathscr{K}$-transform, and the Widder potential transform will be required in
our investigation.
\begin{lem}\label{l1} 
The identities 
\begin{equation}
\label{l1:1}
\mathscr{H}_{\nu}\Big\{\mathscr{K}_{\nu}\big\{f(x);u\big\};y\Big\}=
y^{\nu+1/2}\,\mathscr{P}\big\{x^{-\nu-1/2}\,f(x);y\big\},
\end{equation}
and
\begin{equation}
\label{l1:2}
\mathscr{K}_{\nu}\Big\{\mathscr{H}_{\nu}\big\{f(x);u\big\};y\Big\}=
y^{-\nu+1/2}\,\mathscr{P}\big\{x^{\nu-1/2}\,f(x);y\big\},
\end{equation}
hold true provided that $\Re(\nu)>-1$ and each member of the assertions \eqref{l1:1} and \eqref{l1:2} exists.
\end{lem}
\begin{pf} 
We only give here the proof of the iteration identity \eqref{l1:1} because the proof of \eqref{l1:2} is similar. Indeed, 
by the definitions \eqref{ht} of the Hankel transform and \eqref{kt} of the $\mathscr{K}$-transform, we have 
\begin{align}
\notag
\mathscr{H}_{\nu}\Big\{\mathscr{K}_{\nu}& \big\{f(x);u \big\};y\Big\}\\&=
\int_{0}^{\infty} (u\,y)^{1/2}\,J_{\nu}(u\,y)\, \Big[\int_{0}^{\infty} (u\,x)^{1/2}\,\text{K}_{\nu}(u\,x)\,f(x)\,dx\Big]\,du
\label{l1:p1} 
\end{align}
Changing the order of integration (which is permissible by absolute convergence of the integrals involved), we find from \eqref{l1:p1} that 
\begin{align}
\notag
\mathscr{H}_{\nu}\Big\{\mathscr{K}_{\nu}&\big\{f(x);u\big\};y\Big\}\\&=
\int_{0}^{\infty} (x\,y)^{1/2}\,f(x)\, \Big[\int_{0}^{\infty} u\,J_{\nu}(u\,y)\,\text{K}_{\nu}(u\,x)\,du\Big]\,dx.
\label{l1:p2} 
\end{align}
Using the known formula \cite[p. 658, Entry 6.521-1]{GR}, the integral on the right-hand side of \eqref{l1:p2} is given by
\begin{equation}
\label{l1:p3}
\int_{0}^{\infty} u\,J_{\nu}(u\,y)\,\text{K}_{\nu}(u\,x)\,du=\frac{y^{\nu}}{x^{\nu}\,(x^{2}+y^{2})}\cdot
\end{equation}
The assertion \eqref{l1:1} follows upon substituting the result \eqref{l1:p3} into \eqref{l1:p2} and then using the definition \eqref{wpt} of the Widder potential transform.
\qed
\end{pf}

Using the known property of the Hankel transform
\begin{equation}
\label{ihh}
\mathscr{H}_{\nu}\Big\{\mathscr{H}_{\nu}\big\{f(x);y\big\};x\Big\}=f(x),
\end{equation}
and our identities \eqref{l1:1} and \eqref{l1:2}, we obtain the iteration identities contained in
\begin{cor}\label{c1:l1} 
Under the assumptions of the Lemma \ref{l1},
the  iteration identities 
\begin{gather}
\label{c1:l1:1}
\mathscr{H}_{\nu}\Big\{u^{\nu+1/2}\,\mathscr{P}\big\{x^{-\nu-1/2}\,f(x);u\big\};y\Big\}
=\mathscr{K}_{\nu} \big\{f(x);y \big\}
\\
\label{c1:l1:2}
\mathscr{P}\Big\{u^{\nu-1/2}\,\mathscr{H}_{\nu} \big\{f(x);u \big\};y\Big\}
=
y^{\nu-1/2}\,\mathscr{K}_{\nu} \big\{f(x);y \big\},
\intertext{and}
\label{c1:l1:3}
\mathscr{H}_{\nu}\Big\{u^{\nu+1/2}\,\mathscr{P}\big\{x^{-\nu-1/2}\,f(x);u\big\};y\Big\}
=
y^{-\nu+1/2}\,\mathscr{P}\Big\{u^{\nu-1/2}\,\mathscr{H}_{\nu} \big\{f(x);u \big\};y\Big\}.
\end{gather}
\end{cor}

Setting $\nu=1/2$ and $\nu=-1/2$ in the results \eqref{l1:1} and \eqref{l1:2}; and using the relationships \eqref{ht:fst}, \eqref{ht:fct} and \eqref{kt:lt},  we obtain the results contained in
\begin{cor}\label{c2:l1} 
Under the assumptions of the Lemma \ref{l1},
the following iteration identities 
\begin{align}
\label{c2:l1:1}
\mathscr{F}_{S}\Big\{\mathscr{L}\big\{f(x);u \big\};y\Big\}&=
y\,\mathscr{P}\big\{x^{-1}\,f(x);y\big\},
\\
\label{c2:l1:2}
\mathscr{L}\Big\{\mathscr{F}_{S} \big\{f(x);u \big\};y\Big\}&=
\mathscr{P}\big\{f(x);y\big\},
\\
\label{c2:l1:3}
\mathscr{F}_{C}\Big\{\mathscr{L}\big\{f(x);u \big\};y\Big\}&=
\mathscr{P}\big\{f(x);y\big\}
\end{align}
and
\begin{equation}
\label{c2:l1:4}
\mathscr{L}\Big\{\mathscr{F}_{C} \big\{f(x);u \big\};y\Big\}=
y\,\mathscr{P}\big\{x^{-1}\,f(x);y\big\}
\end{equation}
hold true.
\end{cor}

\begin{rem}\label{r1:l1}
The identity \eqref{c2:l1:2} was obtained earlier in Widder \cite{Wd71}. Thus, our result \eqref{l1:2} in the Lemma \ref{l1} generalizes the earlier result \eqref{c2:l1:2}.
\end{rem}  

Immediate consequences of the Corollary \ref{c2:l1} are contained in 
\begin{cor}\label{c3:l1} 
Under the assumptions of the Lemma \ref{l1},
the following iteration identities 
\begin{equation}
\label{c3:l1:1}
\mathscr{F}_{S}\Big\{\mathscr{L}\big\{f(x);u \big\};y\Big\}=
\mathscr{L}\Big\{\mathscr{F}_{C} \big\{f(x);u \big\};y\Big\}
\end{equation}
and
\begin{equation}
\label{c3:l1:2}
\mathscr{L}\Big\{\mathscr{F}_{S} \big\{f(x);u \big\};y\Big\}=
\mathscr{F}_{C}\Big\{\mathscr{L}\big\{f(x);u \big\};y\Big\}
\end{equation}
hold true.
\end{cor}

Similary, setting $\nu=1/2$ and $\nu=-1/2$ in the results \eqref{c1:l1:1} and \eqref{c1:l1:2}; and using the relationships \eqref{ht:fst}, \eqref{ht:fct} and \eqref{kt:lt},   we obtain the results contained in the following two corollaries.
\begin{cor}\label{c4:l1} 
Under the assumptions of the Lemma \ref{l1},
the following iteration identities 
\begin{gather}
\label{c4:l1:1}
\mathscr{F}_{S}\Big\{u\,\mathscr{P}\big\{x^{-1}\,f(x);u \big\};y\Big\}=
\frac{\pi}{2}\,\mathscr{L}\big\{f(x);y\big\},
\\
\label{c4:l1:2}
\mathscr{P}\Big\{\mathscr{F}_{S} \big\{f(x);u \big\};y\Big\}=
\frac{\pi}{2}\,\mathscr{L}\big\{f(x);y\big\},
\\
\label{c4:l1:3}
\mathscr{F}_{C}\Big\{\mathscr{P}\big\{f(x);u \big\};y\Big\}=
\frac{\pi}{2}\,\mathscr{L}\big\{f(x);y\big\}
\intertext{and}
\label{c4:l1:4}
y\,\mathscr{P}\Big\{u^{-1}\,\mathscr{F}_{C} \big\{f(x);u \big\};y\Big\}=
\frac{\pi}{2}\,\mathscr{L}\big\{f(x);y\big\}
\end{gather}
hold true.
\end{cor}

Immediate consequences of the Corollary \ref{c4:l1} are contained in 
\begin{cor}\label{c5:l1} 
Under the assumptions of the Lemma \ref{l1},
the following iteration identities 
\begin{align}
\notag
\mathscr{F}_{S}\Big\{u\,\mathscr{P}\big\{x^{-1}\,f(x);u \big\};y\Big\}
&=
\mathscr{P}\Big\{\mathscr{F}_{S} \big\{f(x);u \big\};y\Big\}
\\
\notag
&=
\mathscr{F}_{C}\Big\{\mathscr{P}\big\{f(x);u \big\};y\Big\}
\\
\label{c5:l1:1}
&=
y\,\mathscr{P}\Big\{u^{-1}\,\mathscr{F}_{C} \big\{f(x);u \big\};y\Big\}
\end{align}
hold true.
\end{cor}

Our main result is contained in
\begin{thm}\label{t1} 
The following Parseval-Goldstein type identities hold true:
\begin{align}
\label{t1:1}
\int_{0}^{\infty}
\mathscr{H}_{\nu}\big\{f(x);y \big\}
\,\mathscr{K}_{\nu}\big\{g(u);y \big\}\,dy
&=
\int_{0}^{\infty}
x^{\nu+1/2}\,f(x)\,\mathscr{P}\big\{u^{-\nu-1/2}\,g(u);x\big\}\,dx,
\\
\label{t1:2}
\int_{0}^{\infty}
\mathscr{H}_{\nu}\big\{f(x);y \big\}
\,\mathscr{K}_{\nu}\big\{g(u);y \big\}\,dy
&=
\int_{0}^{\infty}
u^{-\nu+1/2}\,g(u)\,\mathscr{P}\big\{x^{\nu-1/2}\,f(x);u\big\}\,du
\intertext{and}
\label{t1:3}
\int_{0}^{\infty}
x^{\nu+1/2}\,f(x)\,\mathscr{P}\big\{u^{-\nu-1/2}\,g(u);x\big\}\,dx
&=
\int_{0}^{\infty}
u^{-\nu+1/2}\,g(u)\,\mathscr{P}\big\{x^{\nu-1/2}\,f(x);u\big\}\,du,
\end{align} 
provided that $\Re(\nu)>-1$ and the integrals involved converge
absolutely.
\end{thm}
\begin{pf} 
We only give here the proof of the Parseval-Goldstein identity \eqref{t1:1} because the proof of \eqref{t1:2} is similar and the identity \eqref{t1:3} follows easily from the assertions \eqref{t1:1} and \eqref{t1:2}. Indeed, by the definition \eqref{ht} of the Hankel transform, we have
\begin{align}
\notag
\int_{0}^{\infty}
\mathscr{H}_{\nu}&\big\{f(x);y \big\}
\,\mathscr{K}_{\nu}\big\{g(u);y \big\}\,dy
\\
\label{t1:p1}&=
\int_{0}^{\infty}
\mathscr{K}_{\nu}\big\{g(u);y \big\}\,
\Big[\int_{0}^{\infty} (x\,y)^{1/2}\,J_{\nu}(x\,y)\,f(x)\,dx\Big]\,dy
\end{align} 
Changing the order of integration which is permissible under the assumptions
of the theorem and using the definition \eqref{ht} of the Hankel transform once again, we find from \eqref{t1:p1} that 
\begin{align}
\notag
\int_{0}^{\infty}
\mathscr{H}_{\nu}&\big\{f(x);y \big\}
\,\mathscr{K}_{\nu}\big\{g(u);y \big\}\,dy
\\
\notag
&=
\int_{0}^{\infty}
f(x)\,\bigg[\int_{0}^{\infty} (x\,y)^{1/2}\,\mathscr{K}_{\nu}\big\{g(u);y \big\}\,
J_{\nu}(x\,y)\,dy\bigg]\,dx
\\
\label{t1:p2}
&=\int_{0}^{\infty}
f(x)\,\mathscr{H}_{\nu}\Big\{\mathscr{K}_{\nu}\big\{g(u);y\big\};x\Big\}\,dx
\end{align} 
Now the assertion \eqref{t1:1} easily follows \eqref{t1:p2} and \eqref{l1:1} of the Lemma \ref{l1}
\qed\end{pf}

Setting $\nu=1/2$ in the Parseval-Goldstein type identities \eqref{t1:1}, \eqref{t1:2}, and \eqref{t1:3} of our Theorem \ref{t1} and using the special cases \eqref{ht:fst} and \eqref{kt:lt} of the Hankel transforms and the $\mathscr{K}$-transforms, respectively; we obtain the identities contained in the following corollary.
\begin{cor}\label{c1:t1} 
The following Parseval-Goldstein type identities hold true:
\begin{align}
\label{c1:t1:1}
\int_{0}^{\infty}
\mathscr{F}_{S}\big\{f(x);y\big\}\,\mathscr{L}\big\{g(u);y\big\}\,dy&=
\int_{0}^{\infty} x\,f(x)\,\mathscr{P}\big\{u^{-1}\,g(u);x\big\}\,dx
\\
\label{c1:t1:2}
\int_{0}^{\infty}
\mathscr{F}_{S}\big\{f(x);y\big\}\,\mathscr{L}\big\{g(u);y\big\}\,dy&=
\int_{0}^{\infty} g(u)\,\mathscr{P}\big\{f(x);u\big\}\,du
\intertext{and}
\label{c1:t1:3}
\int_{0}^{\infty} x\,f(x)\,\mathscr{P}\big\{u^{-1}\,g(u);x\big\}\,dx&
=
\int_{0}^{\infty} g(u)\,\mathscr{P}\big\{f(x);u\big\}\,du
\end{align} 
provided that integrals involved converge absolutely.
\end{cor}
\begin{rem}
\label{r1:t1}
The Parseval-Goldstein type identity \eqref{c1:t1:2} is obtained earlier in Srivastava and Y\"urekli
\cite[p. 586, Eq. (7)]{SY91}. Therefore, the identity \eqref{t1:2} is a
generalization of the earlier identity \eqref{c1:t1:2}.
\end{rem}

Setting $\nu=-1/2$ in the Parseval-Goldstein type identities \eqref{t1:1}, \eqref{t1:2}, and \eqref{t1:3} of our Theorem \ref{t1} and using the special cases \eqref{ht:fct} and \eqref{kt:lt} of the Hankel transforms and the $\mathscr{K}$-transforms, respectively; we obtain the identities contained in the following corollary.
\begin{cor}\label{c2:t1} 
The following Parseval-Goldstein type identities hold true:
\begin{align}
\label{c2:t1:1}
\int_{0}^{\infty}
\mathscr{F}_{C}\big\{f(x);y\big\}\,\mathscr{L}\big\{g(u);y\big\}\,dy&=
\int_{0}^{\infty} f(x)\,\mathscr{P}\big\{g(u);x\big\}\,dx
\intertext{and}
\label{c2:t1:2}
\int_{0}^{\infty}
\mathscr{F}_{C}\big\{f(x);y\big\}\,\mathscr{L}\big\{g(u);y\big\}\,dy&=
\int_{0}^{\infty} u\,g(u)\,\mathscr{P}\big\{x^{-1}\,f(x);u\big\}\,du
\end{align} 
provided that integrals involved converge absolutely.
\end{cor}

\begin{rem}
\label{r2:t1}
If we use the Parseval-Goldstein type identities \eqref{c2:t1:1} and \eqref{c2:t1:2} of Corollary \ref{c2:t1}, we obtain a Goldstein type relationship identical to \eqref{c1:t1:3} of Corollary \eqref{c1:t1}.
\end{rem}

Either applying the known property \eqref{ihh} of the Hankel transform in the Parseval-Goldstein type identities \eqref{t1:1} and \eqref{t1:2} of our Theorem \ref{t1} or directly using the iteration identities \eqref{c1:l1:1} and \eqref{c1:l1:2}, we obtain a new set of Parseval-Goldstein type identities contained in
\begin{thm}\label{t2} 
The following Parseval-Goldstein type identities hold true:
\begin{gather}
\label{t2:1}
\int_{0}^{\infty}
y^{\nu+1/2}\,\mathscr{H}_{\nu}\big\{f(x);y \big\}
\,\mathscr{P}\big\{u^{-\nu-\frac{1}{2}}\,g(u);y \big\}\,dy
=
\int_{0}^{\infty}
f(x)\,\mathscr{K}_{\nu}\big\{g(u);x \big\}\,dx,
\\
\label{t2:2}
\int_{0}^{\infty}
y^{\nu+1/2}\,\mathscr{H}_{\nu}\big\{f(x);y \big\}
\,\mathscr{P}\big\{u^{-\nu-\frac{1}{2}}\,g(u);y \big\}\,dy
=
\int_{0}^{\infty}
g(u)\,\mathscr{K}_{\nu}\big\{f(x);u \big\}\,du,
\intertext{and}
\label{t2:3}
\int_{0}^{\infty}
f(x)\,\mathscr{K}_{\nu}\big\{g(u);x \big\}\,dx
=
\int_{0}^{\infty}
g(u)\,\mathscr{K}_{\nu}\big\{f(x);u \big\}\,du
\end{gather} 
provided that $\Re(\nu)>-1$ and the integrals involved converge
absolutely.
\end{thm}
Setting  $\nu=1/2$ and $\nu=-1/2$ in the Parseval-Goldstein type identities \eqref{t2:1}, \eqref{t2:2}, and \eqref{t2:3} of our Theorem \ref{t2} and using the special cases \eqref{ht:fst}, \eqref{ht:fct} and \eqref{kt:lt} of the Bessel functions, we obtain the identities contained in the following corollary.
\begin{cor}\label{c1:t2} 
The following Parseval-Goldstein type identities hold true:
\begin{gather}
\label{c1:t2:1}
\int_{0}^{\infty}
y\,\mathscr{F}_{S}\big\{f(x);y\big\}\,\mathscr{P}\big\{u^{-1}\,g(u);y\big\}\,dy
=
\frac{\pi}{2}\,\int_{0}^{\infty} f(x)\,\mathscr{L}\big\{g(u);x\big\}\,dx,
\\
\label{c1:t2:2}
\int_{0}^{\infty}
y\,\mathscr{F}_{S}\big\{f(x);y\big\}\,\mathscr{P}\big\{u^{-1}\,g(u);y\big\}\,dy
=
\frac{\pi}{2}\,\int_{0}^{\infty} g(u)\,\mathscr{L}\big\{f(x);u\big\}\,du,
\\
\label{c1:t2:3}
\int_{0}^{\infty}
\mathscr{F}_{C}\big\{f(x);y\big\}\,\mathscr{P}\big\{g(u);y\big\}\,dy
=
\frac{\pi}{2}\,\int_{0}^{\infty} f(x)\,\mathscr{L}\big\{g(u);x\big\}\,dx,
\\
\label{c1:t2:4}
\int_{0}^{\infty}
\mathscr{F}_{C}\big\{f(x);y\big\}\,\mathscr{P}\big\{g(u);y\big\}\,dy
=
\frac{\pi}{2}\,\int_{0}^{\infty} g(u)\,\mathscr{L}\big\{f(x);u\big\}\,du,
\intertext{and}
\label{c1:t2:5}
\int_{0}^{\infty} f(x)\,\mathscr{L}\big\{g(u);x\big\}\,dx
=
\int_{0}^{\infty} g(u)\,\mathscr{L}\big\{f(x);u\big\}\,du,
\end{gather} 
provided that integrals involved converge absolutely.
\end{cor}
\section{A set of useful corollaries}

Several interesting consequences of the results in the previous section will be presented in this section. 
\begin{cor}\label{c2:t2} 
The iteration identities hold true for the Widder potential transform, the $\mathscr{K}_{\nu}$-transform:
\begin{equation}
\label{c2:t2:1}
\mathscr{P}\Big\{x^{2\nu}\,\mathscr{P}\big\{u^{-\nu-\frac{1}{2}}\,g(u);x\big\};t\Big\}
=t^{\nu-\frac{1}{2}}\,\mathscr{K}_{\nu}\Big\{\mathscr{K}_{\nu}\big\{g(u);y\big\};t\Big\},
\end{equation} 
where $-1<\Re(\nu)<3/2$; and for the Widder potential transform and the classical Laplace tansform
\begin{align}
\label{c2:t2:2}
\mathscr{P}\Big\{x\,\mathscr{P}\big\{u^{-1}\,g(u);x\big\};t\Big\}
&=\frac{\pi}{2}\,\mathscr{L}\Big\{\mathscr{L}\big\{g(u);y\big\};t\Big\},
\\
\label{c2:t2:3}
\mathscr{P}\Big\{x^{-1}\,\mathscr{P}\big\{g(u);x\big\};t\Big\}
&=\frac{\pi}{2}\,t^{-1}\,\mathscr{L}\Big\{\mathscr{L}\big\{g(u);y\big\};t\Big\}, 
\end{align}
provided that each member of the assertion \eqref{c2:t2:1}, \eqref{c2:t2:2}, and \eqref{c2:t2:3} exists.
\end{cor}
 
\begin{pf}
We set 
\begin{equation} 
\label{c2:t2:p1}
f(x)=x^{\nu+\frac{1}{2}}\,(x^{2}+t^{2})^{-1}
\end{equation}
in the Parseval-Goldstein identity \eqref{t1:1} of our Theorem \ref{t1}. Using
the known formula \cite[p. 35, Entry 4.13]{Y89} yields 
\begin{equation}
\label{c2:t2:p2}
\mathscr{H}_{\nu}\big\{x^{\nu+\frac{1}{2}}\,(x^{2}+t^{2})^{-1};y\big\}=t^{\nu}\,y^{1/2}\,K_{\nu}(t\,y),
\end{equation}
where $-1<\Re(\nu)<3/2$.
The assertion \eqref{c2:t2:1} follows upon substituting \eqref{c2:t2:p1} and \eqref{c2:t2:p2} into the identity \eqref{t1:1} and using the definitions \eqref{kt} and \eqref{wpt} of the Widder potential transform and the $\mathscr{K}_{\nu}$-transform, respectively. The assertions \eqref{c2:t2:2} and \eqref{c2:t2:3} immediately follows when we set $\nu=1/2$ and $\nu=-1/2$ in \eqref{c2:t2:1}, respectively; and using the special cases \eqref{kt:lt}.
\end{pf}

\begin{rem}\label{r2:t2} 
It is well known that the second iterate of the classical Laplace transform is the Stieltjes transrform; that is,
\begin{equation}
\label{r2:t2:1}
\mathscr{L}\Big\{\mathscr{L}\big\{g(u);x\big\};t\Big\}
=\mathscr{S}\big\{g(u);y\big\}.
\end{equation} 
Using the result \eqref{r2:t2:1} in the identities \eqref{c2:t2:2} and \eqref{c2:t2:3} of our Corollary \ref{c2:t2}, we deduce
\begin{align}
\label{r2:t2:2}
\mathscr{P}\Big\{x\,\mathscr{P}\big\{u^{-1}\,g(u);x\big\};t\Big\}
&=\frac{\pi}{2}\,\mathscr{S}\big\{g(u);y\big\}, 
\\
\label{r2:t2:3}
\mathscr{P}\Big\{x^{-1}\,\mathscr{P}\big\{g(u);x\big\};t\Big\}
&=\frac{\pi}{2}\,t^{-1}\,\mathscr{S}\big\{g(u);y\big\}, 
\end{align}
provided that each member of the assertion \eqref{r2:t2:2} and \eqref{r2:t2:3} exists.
\end{rem}

\begin{cor}\label{c3:t2} 
If $-\Re(\nu)-3/2<\Re(\mu)<-1/2$, then the following Parseval type identities hold true for the Widder transform and  the $\mathscr{K}$-transform:
\begin{align}
\notag
\int_{0}^{\infty}x^{\mu+\nu+\frac{1}{2}}\,&\mathscr{P}\big\{u^{-\nu-\frac{1}{2}}\,g(u);x\big\}\,dx
\\
\label{c3:t2:1}
&=
2^{\mu+\frac{1}{2}}\,
\frac{\Gamma\big(\frac{\mu}{2}+ \frac{\nu}{2}+ \frac{3}{4}\big)}
{\Gamma\big(\frac{\nu}{2}-\frac{\mu}{2}+ \frac{1}{4}\big)}\,
\int_{0}^{\infty}y^{-\mu-1}\,\mathscr{K}_{\nu}\big\{g(u);y\big\}\,dy,
\\
\notag
\int_{0}^{\infty}y^{-\mu-1}&\,\mathscr{K}_{\nu}\big\{g(u);y\big\}\,dy
\\
\label{c3:t2:2}&=
2^{-\mu-\frac{3}{2}}\,
\Gamma\Big(\frac{\nu}{2}-\frac{\mu}{2}+\frac{1}{4}\Big)\,\Gamma\Big(\frac{1}{4}-\frac{\mu}{2}-\frac{\nu}{2}\Big)
\,
\int_{0}^{\infty}u^{\mu}\,g(u)\,du
\intertext{and}
\notag
\int_{0}^{\infty}x^{\mu+\nu+\frac{1}{2}}&\,\mathscr{P}\big\{u^{-\nu-\frac{1}{2}}\,g(u);x\big\}\,dx
\\
\label{c3:t2:3}
&=
\frac{1}{2}\,
\Gamma\big(\frac{\mu}{2}+\frac{\nu}{2}+\frac{3}{4}\big)\,\Gamma\Big(\frac{1}{4}-\frac{\mu}{2}-\frac{\nu}{2}\Big)
\,
\int_{0}^{\infty}u^{\mu}\,g(u)\,du
\end{align} 
provided that each member of the assertion \eqref{c2:t2:1}, \eqref{c2:t2:2}, and \eqref{c2:t2:3} exists.
\end{cor}
 
\begin{pf}
We set 
\begin{equation} 
\label{c3:t2:p1}
f(x)=x^{\mu}
\end{equation}
in the Theorem \ref{t1}. Using the known formula \cite[p. 248, Entry (A1)]{Y89} yields
\begin{equation}
\label{c3:t2:p2}
\mathscr{P}\big\{x^{\mu+\nu-\frac{1}{2}};u\big\}=\frac{\pi}{2}\,\sec\Big[\frac{\pi}{2}\Big(\mu+\nu+\frac{1}{2}\Big)\Big]\,u^{\mu+\nu-\frac{1}{2}}.
\end{equation}
 It is known that (cf. \cite[p. 22, Entry (7)]{E2}) 
\begin{equation}
\label{c3:t2:p3}
\mathscr{H}_{\nu}\big\{x^{\mu};y\big\}=2^{\mu+\frac{1}{2}}\,y^{-\mu-1}\,
\frac{\Gamma\big(\frac{\mu}{2}+\frac{\nu}{2}+\frac{3}{4}\big)}
{\Gamma\big(\frac{\nu}{2}-\frac{\mu}{2}+\frac{1}{4}\big)}\,\cdot
\end{equation}
Now the assertions \eqref{c3:t2:1} and \eqref{c3:t2:2} follows upon substituting \eqref{c3:t2:p1}, \eqref{c3:t2:p2} and \eqref{c3:t2:p3} into the Parseval-Goldstein identities \eqref{t1:1} and \eqref{t1:2}, and then utilizing the formula for the gamma function
\begin{equation}
\label{c3:t2:p4}
\Gamma\Big(\frac{1}{2}+z\Big)\,\Gamma\Big(\frac{1}{2}-z\Big)=
\frac{\pi}{\cos(\pi\,z)}\,\cdot
\end{equation} 
The assertion \eqref{c3:t2:1} immediately follows from \eqref{c3:t2:1} and \eqref{c3:t2:2}.
\end{pf}

Setting $\nu=1/2$ in the Corollary \ref{c3:t2} and utilizing the formulas \eqref{c3:t2:p4} and 
\begin{equation}
\label{gf:1}
\Gamma(2\lambda)=\sqrt{\pi}\,2^{2\lambda-1}\,\Gamma(\lambda)\,\Gamma\Big(\lambda+\frac{1}{2}\Big)
\end{equation}
for the Gamma function, we obtain the following Parseval type identities for the Laplace transforms as special cases of the identities given in the previous corollary.

\begin{cor}\label{c4:t2} 
Each of the following identities holds true:
\begin{align}
\notag
\int_{0}^{\infty}y^{-\mu-1}&\,\mathscr{L}\big\{g(u);y \big\}\,dy
\\
\label{c4:t2:1}
&=
\sqrt{\frac{\pi}{2}}\,\sec\Big(\frac{\pi\,\mu}{2}\Big)\,\frac{1}{\Gamma(\mu+1)}\,\int_{0}^{\infty}x^{\mu+1}\,\mathscr{P} \big\{u^{-1}\,g(u);x\big\}\,dx
\\
\notag
&(\Re(\mu)>-1),\\
\label{c4:t2:2}
\int_{0}^{\infty}y^{-\mu-1}&\,\mathscr{L} \big\{g(u);y \big\}\,dy
=
\Gamma(-\mu)\,\int_{0}^{\infty}u^{\mu}\,g(u)\,du
\qquad (\Re(\mu)<0),
\intertext{and}
\label{c4:t2:3}
\int_{0}^{\infty}x^{\mu+1}&\,\mathscr{P} \big\{u^{-1}\,g(u);x\big\}\,dx=
\sqrt{\frac{2}{\pi}}\,\cos\Big(\frac{\pi\,\mu}{2}\Big)\,\frac{\Gamma(-\mu)}{\Gamma(\mu+1)}\,
\int_{0}^{\infty}u^{\mu}\,g(u)\,du
\\
\notag
&(-1<\Re(\mu)<0),
\end{align} 
provided that each member of the assertion \eqref{c4:t2:1}, \eqref{c4:t2:2}, and \eqref{c4:t2:3} exists.
\end{cor}
\begin{rem}\label{r3:t2} 
Setting $\mu=-1$ in the identity \eqref{c4:t2:2} of our Corollary \ref{c4:t2} we obtain the well known identity for the Laplace transform (cf. \cite[p. 110, Eq. (2.4)]{HMRY})
\begin{equation}
\label{r3:t2:1}
\int_{0}^{\infty} \mathscr{L}\big\{g(u);u\big\}\,du=\int_{0}^{\infty}\frac{g(u)}{u}\,du.
\end{equation}
Hence the identity \eqref{c2:t2:2} is a generalization of the known identity \eqref{r3:t2:1}.
\end{rem}

\section{Illustrative Examples}

Several interesting consequences of the results in the previous section will be presented in this section. 
An interesting illustration for the identity \eqref{l1:2} of our Lemma \ref{l1} is contained in the following example.
\begin{exmp}\label{e1}
 We show that
\begin{equation}
\label{e1:1}
\mathscr{P}\bigg\{\frac{x^{2\nu}}{x^{2}+a^{2}};y \bigg\}=\frac{\pi}{2\,\sin(\nu\,\pi)}\,\frac{a^{2\nu}-y^{2\nu}}{a^{2}-y^{2}}, \qquad \big(\Re(y+a)>0, \quad |\Re(\nu)|<1\big).
\end{equation}
\end{exmp}
\begin{pf}
 We put 
\begin{equation}
\label{e1:p1}
f(x)=\frac{x^{\nu+\frac{1}{2}}}{x^{2}+a^{2}}\cdot
\end{equation}
Utilizing the known formula \cite[p. 23, Entry (12)]{E2}, we have
\begin{equation}
\label{e1:p2}
\mathscr{H}_{\nu}\bigg\{\frac{x^{\nu+\frac{1}{2}}}{x^{2}+a^{2}};u\bigg\}=a^{\nu}\,u^{1/2}\,\text{K}_{\nu}(a\,u).
\end{equation}
Applying the $\mathscr{K}$-transform to both sides of the equation \eqref{e1:p2} and then using the known formula  \cite[p. 145, Entry (48)]{E2}
\begin{align}
\notag
\mathscr{K}_{\nu}\Bigg\{\mathscr{H}_{\nu}\bigg\{\frac{x^{\nu+\frac{1}{2}}}{x^{2}+a^{2}};u\bigg\};y\Bigg\}&=a^{\nu}\,\mathscr{K}_{\nu}\big\{u^{1/2}\,\text{K}_{\nu}(a\,u);y\big\}
\\
\label{e1:p3}
&=\frac{\pi\,y^{\frac{1}{2}-\nu}}{2\,\sin(\nu\,\pi)}\,\frac{a^{2\nu}-y^{2\nu}}{a^{2}-y^{2}}\cdot
\end{align}
The assertion \eqref{e1:1} follows upon substituting \eqref{e1:p1} and  \eqref{e1:p3} into the iteration identity \eqref{l1:2} of our Lemma \ref{l1}.
\qed\end{pf}

\begin{rem}\label{r1:e1} 
Using the definition \eqref{wpt} of the Widder transform, we can restate the result \eqref{e1:1} as \begin{equation}
\label{r1:e1:1}
\int_{0}^{\infty} \frac{x^{2\nu+1}}{\big(x^{2}+a^{2}\big)\big(x^{2}+y^{2}\big)}\,dx=
\frac{\pi}{2\,\sin(\nu\,\pi)}\,\frac{a^{2\nu}-y^{2\nu}}{a^{2}-y^{2}}.
\end{equation}
Setting $\beta=a^{2}$, $\gamma=y^{2}$, and $\nu=\frac{\mu}{2}-1$ we obtain the known formula
\begin{align}
\notag
\int_{0}^{\infty} \frac{x^{\mu-1}}{\big(x^{2}+\beta\big)\big(x^{2}+\gamma\big)}\,dx=
\frac{\pi}{2\,\sin(\mu\,\pi)}\,\frac{\gamma^{\frac{\mu}{2}-1}-\beta^{\frac{\mu}{2}-1}}{\beta-\gamma},
\\
\label{r1:e1:2}
\big(|\text{\rm arg}(\beta)|<\pi,\quad |\text{\rm arg}(\gamma)|<\pi,\quad 0<\Re(\mu)<4\big),
\end{align}
(cf. \cite[p. 300, Entry (3.264-2)]{GR}).

The Mellin transform is defined as 
\begin{equation}
\label{mt}
\mathscr{M}\big\{f(x);\mu\big\}=\int_{0}^{\infty}x^{\mu-1}\,f(x)\,dx.
\end{equation}
Using the formula \eqref{r1:e1:2} and the definition \eqref{mt} of the Mellin transform, we deduce another known formula \cite[p. 309, Entry (14)]{E1}. These results verify our results in the Lemma \ref{l1}.
\end{rem}

\begin{exmp}\label{e2}
 We show that
\begin{align}
\notag
\int_{0}^{\infty}x^{\mu+\nu+\frac{1}{2}}\,&\frac{a^{2\nu}-x^{2\nu}}{a^{2}-x^{2}}\,dx
=\frac{1}{2}\,\sin(\nu\,\pi)\,
\\
\notag
&\sec\Big[\pi\Big(\frac{\mu}{2}+\frac{3\nu}{2}+\frac{1}{4}\Big)\Big]
\,\Gamma\Big(\frac{\mu}{2}+\frac{\nu}{2}+\frac{3}{4}\Big)
\,\Gamma\Big(\frac{1}{4}-\frac{\mu}{2}-\frac{\nu}{2}\Big)\,a^{\mu+3\nu-\frac{1}{2}}
\\
\label{e2:1}
&\qquad \big(\Re(y+a)>0, \quad |\Re(\nu)|<1,\quad -\frac{3}{2}<\Re(\mu+3\nu)<\frac{1}{2}\big).
\end{align}
\end{exmp}
\begin{pf}
 We put 
\begin{equation}
\label{e2:p1}
g(u)=\frac{u^{3\nu+\frac{1}{2}}}{u^{2}+a^{2}}\cdot
\end{equation}
Using the formula \eqref{e1:1} of our Example \ref{e1}, we have
\begin{equation}
\label{e2:p2}
\mathscr{P}\bigg\{u^{-\nu-\frac{1}{2}}\,g(u);x\bigg\}=
\mathscr{P}\bigg\{\frac{u^{2\nu}}{u^{2}+a^{2}};x\bigg\}
=
\frac{\pi}{2\sin(\nu\,\pi)}\,\frac{a^{2\nu}-x^{2\nu}}{a^{2}-x^{2}}.
\end{equation}
Using the known formula \cite[p. 248, Entry (A1)]{Y89} yields
\begin{equation}
\label{e2:p3}
\mathscr{P}\big\{x^{\mu+3\nu-\frac{1}{2}};u\big\}=\frac{\pi}{2}\,\sec\Big[\frac{\pi}{2}\Big(\mu+3\nu+\frac{1}{2}\Big)\Big]\,u^{\mu+3\nu-\frac{1}{2}}.
\end{equation}
We obtain the assertion \eqref{e2:1} upon substituting the equations \eqref{e2:p1}, \eqref{e2:p2} and \eqref{e2:p3} into the identity \eqref{c3:t2:3} of our Corollary \eqref{c3:t2}.
 \qed\end{pf}
\begin{rem}\label{r1:e2} 
Setting $\nu=1/2$ in the formula \eqref{e2:1} and using the formulas
\begin{equation}
\label{r1:e2:1}
\Gamma(z+1)=z\,\Gamma(z),\qquad \Gamma(z)\,\Gamma(-z)=\frac{-\pi}{z\,\sin(\pi\,z)}
\end{equation}
for the gamma function, we obtain the known integral formula
\begin{equation}
\label{r1:e2:2}
\int_{0}^{\infty} \frac{x^{\mu+1}}{x+a}\,dx=\pi\,\csc(\pi\,\mu)\,a^{\mu+1}\qquad 
\big(a>0, \quad -2<\Re(\mu)<-1\big),
\end{equation}
(cf. \cite[p. 289, Entry 3.222-2]{GR}). Thus, our formula \eqref{e2:1} is a generalization of the known formula \eqref{r1:e2:2}.
\end{rem}

The following example contains a result involving Struve's functions $\mathbf{H}_{\nu}(x)$ (see \cite[p. 372]{E1}).
\begin{exmp}\label{e3}
 We show that
\begin{align}
\notag
\int_{0}^{\infty}u^{\mu+\frac{1}{2}}\,\mathbf{H}_{\nu}(a\,u)\,du
=\frac{2^{\mu+\frac{1}{2}}\,\Gamma\big(\frac{\mu}{2}+\frac{\nu}{2}+\frac{3}{4}\big)}
{a^{\mu+\frac{3}{2}}\,\Gamma\big(\frac{\nu}{2}-\frac{\mu}{2}+\frac{1}{4}\big)}
\,\tan\Big[\pi\,\Big(\frac{\mu}{2}+\frac{\nu}{2}+\frac{3}{4}\Big)\Big],
\\
\label{e3:1}
\qquad \big(\Re(\nu)>-\frac{3}{2},\quad -\frac{5}{2}<\Re(\mu+\nu)<-\frac{1}{2}\big).
\end{align}
\end{exmp}
\begin{pf}
 We put 
\begin{equation}
\label{e3:p1}
g(u)=u^{\frac{1}{2}}\,\mathbf{H}_{\nu}(a\,u)
\end{equation}
Using the known formula \cite[p. 150, Entry (73)]{E2}, we have
\begin{equation}
\label{e3:p2}
\mathscr{K}_{\nu}\big\{u^{\frac{1}{2}}\,\mathbf{H}_{\nu}(a\,u);x\big\}=
\frac{a^{\nu+1}\,y^{-\nu-\frac{1}{2}}}{y^{2}+a^{2}}\cdot
\end{equation}
Utilizing the formula \cite[p. 248, Entry (A1)]{Y89} yields
\begin{equation}
\label{e3:p3}
\mathscr{P}\big\{y^{-\mu-\nu-\frac{5}{2}};a\big\}=\frac{\pi}{2}\,\sec\Big[\frac{\pi}{2}\Big(\mu+\nu+\frac{3}{2}\Big)\Big]\,a^{-\mu-\nu-\frac{5}{2}}.
\end{equation}
Substituting the results \eqref{e3:p1}, \eqref{e3:p2} and \eqref{e3:p3} into the identity \eqref{c3:t2:2} of our Corollary \eqref{c3:t2}, we obtain
\begin{equation}
\label{e3:p4}
\int_{0}^{\infty}u^{\mu+\frac{1}{2}}\,\mathbf{H}_{\nu}(a\,u)\,du=
\frac{2^{\mu+\frac{1}{2}}\,\pi\,\sec\big[\frac{\pi}{2}\big(\mu+\nu+\frac{3}{2}\big)\big]}
{a^{\mu+\frac{3}{2}}\,\Gamma\big(\frac{\nu}{2}-\frac{\mu}{2}+\frac{1}{4}\big)
\,\Gamma\big(\frac{1}{4}-\frac{\mu}{2}-\frac{\nu}{2}\big)}
\cdot
\end{equation}
Multiplying and dividing the fraction on the right hand side of \eqref{e3:p4} with $\Gamma\big(\frac{\mu}{2}+\frac{\nu}{2}+\frac{3}{4}\big)$ and using the formula 
\begin{equation}
\Gamma(z)\,\Gamma(1-z)=\frac{\pi}{\sin(\pi z)}
\end{equation}
with $z=\frac{\mu}{2}+\frac{\nu}{2}+\frac{3}{4}$, we obtain the desired formula \eqref{e3:1}.
 \qed\end{pf}
 
\begin{rem}\label{r1:e3} 
Using the definition \eqref{mt} of the Mellin transform together with \eqref{e3:1} of our Example \eqref{e3}, we obtain 
\begin{equation}
\label{r1:e3:1}
\mathscr{M} \bigg\{\mathbf{H}_{\nu}(a\,u);\mu+\frac{3}{2}\bigg\}
=\frac{2^{\mu+\frac{1}{2}}\,\Gamma\big(\frac{\mu}{2}+\frac{\nu}{2}+\frac{3}{4}\big)}
{a^{\mu+\frac{3}{2}}\,\Gamma\big(\frac{\nu}{2}-\frac{\mu}{2}+\frac{1}{4}\big)}
\,\tan\Big[\pi\,\Big(\frac{\mu}{2}+\frac{\nu}{2}+\frac{3}{4}\Big)\Big],
\end{equation}
(cf. \cite[p. 335, Entry (52)]{E1}). The $\mathbf{H}$-transform is defined as.
\begin{equation}
\label{bfht}
\frak{H}_{\nu}\big\{f(x);y\big\}=\int_{0}^{\infty} \sqrt{x\,y}\, \mathbf{H}_{\nu}(x\,y)\,f(x)\,dx\end{equation}
Using the definition \eqref{bfht} of the $\mathbf{H}$-transform together with \eqref{e3:1}, we deduce
\begin{equation}
\label{r1:e3:2}
\frak{H}_{\nu}\big\{u^{\mu};a\big\}=
\frac{2^{\mu+\frac{1}{2}}\,\Gamma\big(\frac{\mu}{2}+\frac{\nu}{2}+\frac{3}{4}\big)}
{a^{\mu+1}\,\Gamma\big(\frac{\nu}{2}-\frac{\mu}{2}+\frac{1}{4}\big)}
\,\tan\Big[\pi\,\Big(\frac{\mu}{2}+\frac{\nu}{2}+\frac{3}{4}\Big)\Big],
\end{equation}
(cf. \cite[p. 158, Entry (4)]{E2}).
\end{rem}

\begin{exmp}\label{e4}
 We show that
\begin{align}
\notag
\int_{0}^{\infty} \frac{y^{\nu-\mu-\frac{1}{2}}}{\big(y^{2}+a^{2}\big)^{2\nu+1}}\,dy
=
\frac{1}{2}\,\frac{a^{-3\nu-\mu-\frac{3}{2}}}{\Gamma(2\nu+1)}
\,\Gamma\Big(\frac{\nu}{2}-\frac{\mu}{2}+\frac{1}{4}\Big)
\,\Gamma\Big(\frac{3\nu}{2}+\frac{\mu}{2}+\frac{3}{4}\Big)
\\
\label{e4:1}
\qquad \big(\Re(\mu+3\nu)>-\frac{3}{2},\quad \Re(\mu-\nu)<\frac{1}{2}\big).
\end{align}
\end{exmp}
\begin{pf}
 We put 
\begin{equation}
\label{e4:p1}
g(u)=u^{2\nu+\frac{1}{2}}\,\text{J}_{\nu}(a\,u)
\end{equation}
Using the known formula \cite[p. 137, Entry (16)]{E2}, we have
\begin{equation}
\label{e4:p2}
\mathscr{K}_{\nu}\big\{u^{2\nu+\frac{1}{2}}\,\text{J}_{\nu}(a\,u);y\big\}=
\frac{2^{2\nu}\,a^{\nu}\,y^{\nu+\frac{1}{2}}\,
\Gamma(2\nu+1)}{\big(y^{2}+a^{2}\big)^{2\nu+1}}\cdot
\end{equation}
Substituting the results \eqref{e4:p1} and \eqref{e4:p2} into the identity \eqref{c3:t2:2} of our Corollary \eqref{c3:t2}, we obtain
\begin{align}
\notag
\int_{0}^{\infty} \frac{y^{\nu-\mu-\frac{1}{2}}}{\big(y^{2}+a^{2}\big)^{2\nu+1}}\,dy
=&
\frac{2^{-2\nu-\mu-\frac{3}{2}}\,a^{-\nu}}{\Gamma(2\nu+1)}
\,\Gamma\Big(\frac{\nu}{2}-\frac{\mu}{2}+\frac{1}{4}\Big)
\\
\label{e4:p3}
¬&\qquad
\,\Gamma\Big(\frac{1}{4}-\frac{\mu}{2}-\frac{\mu}{2}\Big)
\int_{0}^{\infty}u^{\mu+2\nu+\frac{1}{2}}\,\text{J}_{\nu}(a\,u)\,du.
\end{align}
The integral on the right hand side of \eqref{e4:p3} can be calculated using the definition \eqref{ht} of the Hankel transform  and using the known formula \cite[p. 22, Entry (7)]{E2}, we obtain
\begin{equation}
\label{e4:p4}
\mathscr{H}_{\nu}\big\{u^{\mu+2\nu};a\big\}=
2^{\mu+2\nu+\frac{1}{2}}\,a^{-\mu-2\nu-1}
\frac{\Gamma\big(\frac{\mu}{2}+\frac{3\nu}{2}+\frac{3}{4}\big)}
{\Gamma\big(\frac{1}{4}-\frac{\nu}{2}-\frac{\mu}{2}\big)}\cdot
\end{equation}
The assertion \eqref{e4:1} follows immediately upon substituting \eqref{e4:p4} into \eqref{e4:p3}.
 \qed\end{pf}

\begin{rem}\label{r1:e4} 
Setting $\nu=0$ in \eqref{e4:1} of our Example \ref{e4}, we obtain
\begin{equation}
\label{e4:1}
\int_{0}^{\infty} \frac{y^{\nu-\mu-\frac{1}{2}}}{y^{2}+a^{2}}\,dy
=
\frac{1}{2}\,a^{-3\nu-\mu-\frac{3}{2}}\,
\,\Gamma\Big(\frac{1}{4}-\frac{\mu}{2}\Big)
\,\Gamma\Big(\frac{\mu}{2}+\frac{3}{4}\Big),
\quad \big(-\frac{3}{2}<\Re(\mu)<\frac{1}{2}\big).
\end{equation}
Using the definition \eqref{wpt} of the Widder transform together with the formula \eqref{c3:t2:p4} for the gamma function and setting $\mu=-\rho-\frac{1}{2}$, we obtain the known formula \cite[p. 248, Entry (A1)]{Y89}.
\end{rem}

\begin{exmp}\label{e5}
 We show that
\begin{align}
\notag
\mathscr{K}_{\nu}\big\{x^{\mu+2\nu};a\big\}=
2^{\mu+2\nu-\frac{1}{2}}\,a^{-2\nu-\mu-\frac{1}{2}}
\,
\Gamma\Big(\frac{\mu}{2}+\frac{\nu}{2}+\frac{3}{4}\Big)\,
\Gamma\Big(\frac{\mu}{2}+\frac{3\nu}{2}+\frac{3}{4}\Big)
\\
\label{e5:1}
\qquad \big(\Re(\mu+2\nu)>|\Re(\nu)|-\frac{1}{2}\big).
\end{align}
\end{exmp}
(cf. \cite[p. 127, Entry 1]{E2}).
\begin{pf}
We put the function $g(u)$ defined in \eqref{e4:p1} into the identity \eqref{c3:t2:1} of our Corollary~ \ref{c3:t2}.
Using the known formula \cite[p. 249, Entry (A8)]{Y89}, we have
\begin{equation}
\label{e5:p2}
\mathscr{P}\big\{u^{-\nu-\frac{1}{2}}\,g(u);y\big\}=\mathscr{P}\big\{u^{\nu}\,\text{J}_{\nu}(a\,u);y\big\}=u^{\nu}\,\text{K}_{\nu}(a\,u).
\end{equation}
Substituting the results \eqref{e4:p2} and \eqref{e5:p2} into the identity \eqref{c3:t2:1}, we obtain
\begin{equation}
\label{e5:p4}
\int_{0}^{\infty} x^{\mu+2\nu}\,\text{K}_{\nu}(a\,x)\,dx
=2^{\mu+2\nu-\frac{1}{2}}\,a^{-2\nu-\mu-\frac{1}{2}}\,
\,\Gamma\Big(\frac{\mu}{2}+\frac{\nu}{2}+\frac{3}{4}\Big)
\,\Gamma\Big(\frac{\mu}{2}+\frac{3\nu}{2}+\frac{3}{4}\Big)
\end{equation}
Using the definition \eqref{kt} on the left hand side of \eqref{e5:p4}, we obtain the assertion \eqref{e5:1}.
 \qed\end{pf}

\section{Conclusions and further directions}
We have obtained several iterations identities for Hankel transforms, $\mathcal{K}$-transforms, and Widder transforms. These iteration identities yield new Parseval-Goldstein identities and Goldstein exchange identities for these and many other integral transforms. We illustrated some use of these identities by evaluating some improper integrals without using complicated techniques. Many other infinite integrals can be evaluated in this manner by applying the lemmas, the theorems and its corollaries considered here.

The methods in this work and many other articles included in the references will yield new identities for the known integral transforms. Using these identities it is possible to compute difficult infinite integrals and integral transforms.




\end{document}